\documentclass[fleqn]{mat01}
\usepackage{times,mathtimy,amssymb}
\begin{document}

\setcounter{page}{39}
\firstpage{39}

\font\xx=msam5 at 10pt
\def\ab{\mbox{\xx{\char'03}}}

\def\definit{\trivlist\item[\hskip\labelsep{DEFINITION.}]}
\def\quest{\trivlist\item[\hskip\labelsep{\it Question.}]}
\def\examm{\trivlist\item[\hskip\labelsep{\it Example.}]}

\newtheorem{theo}{Theorem}
\renewcommand\thetheo{\arabic{section}.\arabic{theo}}
\newtheorem{theor}[theo]{\bf Theorem}
\newtheorem{lem}[theo]{Lemma}
\newtheorem{pot}[theo]{Proof of Theorem}
\newtheorem{propo}[theo]{\rm PROPOSITION}
\newtheorem{rema}[theo]{Remark}
\newtheorem{defn}[theo]{\rm DEFINITION}
\newtheorem{exam}{Example}
\newtheorem{coro}[theo]{\rm COROLLARY}
\newtheorem{step}{\it Step}

\title{Height in splittings of hyperbolic groups}

\markboth{Mahan Mitra}{Height in splittings of hyperbolic groups}

\author{MAHAN MITRA}

\address{29F, Creek Row, Kolkata 700~014, India}

\volume{114}

\mon{February}

\parts{1}

\Date{MS received 8 November 2003}

\begin{abstract} 
Suppose $H$ is a hyperbolic subgroup of a hyperbolic group $G$. Assume
there exists $n > 0$ such that the intersection of $n$ essentially
distinct conjugates of $H$ is always finite. Further assume $G$ splits
over $H$ with hyperbolic vertex and edge groups and the two inclusions
of $H$ are quasi-isometric embeddings. Then $H$ is quasiconvex in $G$.
This answers a question of Swarup and provides a partial converse to the
main theorem of \cite{GMRS}.
\end{abstract}

\keyword{Hyperbolic groups; quasi-isometric embeddings; splittings of
groups.}

\maketitle

\renewcommand\thefootnote{}
\footnote{\!\!\!{\bf Editors' comments}\vspace{.4pc}

We thank the referee for the following comments:\vspace{.2pc}

This paper is an unedited publication of Mahan Mitra's 1997 preprint.
The paper has been used or referred to in the following papers and
book:\vspace{.2pc}

G~A~Swarup, Proof of a weak hyperbolization theorem, {\it Quart. J.
Math.} {\bf 51} (2000) 529--533\vspace{.1pc}

M~Kapovich and B~Kleiner, Hyperbolic groups with low-dimensional
boundary, {\it Ann. Sci. de ENS Paris}, t. {\bf 33} (2000) 647--669\vspace{.2pc}

and the Monograph:\vspace{.2pc}

M~Kapovich, Hyperbolic manifolds and discrete groups, Volume~183 in the
series Progress in Mathematics (Boston, Basel, London: Burkhauser)
(2001).\vspace{.2pc}

Ilya Kapovich has published a (slightly weaker) version of the Main
Theorem~4.6 of Mitra's paper in ``The combination theorem and
quasiconvexity'', {\it Int. J. Algebra Comput.} {\bf 11} (2001) 185--216.
The stronger form is used in Swarup's paper referred to above and is
concerned with the notion of acylindricity or `no long annuli'. In \S5
of Mitra's paper, he mentions that his main argument does not generalize
to graphs of groups. But as observed in the above papers, the extension
to graphs of groups is straightforward. It is not clear what Mitra had
in mind. Since the subsequent questions he raises are still interesting
and unsolved, the presentation has not been\break changed.\vspace{.4pc}

\hfill {\it Editors}
}

\section{Introduction}

Let $G$ be a hyperbolic group in the sense of Gromov \cite{Gromov}. Let
$H$ be a hyperbolic subgroup of $G$. We choose a finite symmetric
generating set for $H$ and extend it to a finite symmetric generating
set for $G$. Let $\Gamma_H$ and $\Gamma_G$ denote the Cayley graphs of
$H$, $G$ respectively with respect to these generating sets.

If $H$ is not quasiconvex in $G$, we would like to understand the group
theoretic (or algebraic) mechanism contributing to the distortion of $H$
in $G$. The first examples of distorted hyperbolic subgroups of
hyperbolic groups were fiber subgroups of fundamental groups of closed
hyperbolic 3-manifolds fibering over the circle. The extrinsic geometry
in this case was studied in detail by Cannon and Thurston \cite{CT} and
later by the author \cite{mitra1,mitra2}. General examples of
normal hyperbolic subgroups of hyperbolic groups have been studied in
\cite{BFH,Mosher2}. A substantially larger class of examples
arise from the combination theorem of Bestvina and Feighn \cite{BF}. In
fact almost all examples of distorted hyperbolic subgroups of hyperbolic
groups use the combination theorem in an essential way (see 
\cite{brady,mitra3} however). It is natural to wonder if there are any other
methods of building distorted hyperbolic subgroups. To get a handle on
this issue one needs the notion of height of a subgroup \cite{GMRS}.

\begin{definit}$\left.\right.$\vspace{.5pc}

\noindent Let $H$ be a subgroup of a group $G$. We say that the
elements $\{g_i |1 \le i \le n\}$ of $G$ are essentially distinct if
$Hg_i \neq Hg_j$ for $i \neq j$. Conjugates of $H$ by essentially
distinct elements are called essentially distinct conjugates.\vspace{.5pc}
\end{definit}

Note that we are abusing notation slightly here, as a conjugate of $H$
by an element belonging to the normalizer of $H$ but not belonging to
$H$ is still essentially distinct from $H$. Thus in this context a
conjugate of $H$ records (implicitly) the conjugating element.

\begin{definit}$\left.\right.$\vspace{.5pc}

\noindent We say that the height of an infinite subgroup $H$ in $G$ is
$n$ if there exists a collection of $n$ essentially distinct conjugates
of $H$ such that the intersection of all the elements of the collection
is infinite and $n$ is maximal possible. We define the height of a
finite subgroup to be $0$.\vspace{.5pc}
\end{definit}

The following question of Swarup \cite{bestvinahp} formulates the
problem we would like to address in this paper:

\begin{quest}
Suppose $H$ is a finitely presented subgroup of a hyperbolic group $G$.
If $H$ has finite height, is $H$ quasiconvex in $G$?
A special case to be considered is when $G$ splits over $H$ and
the inclusions are quasi-isometric embeddings.\vspace{.5pc}
\end{quest}

We shall answer the above question affirmatively in the special
case mentioned.

\setcounter{section}{4}
\setcounter{theo}{5}
\begin{theor}[\!]
Let $G$ be a hyperbolic group splitting over $H$ {\rm (}i.e. $G =
{G_1}{*_H}{G_2}$ or $G = {G_1}{*_H}${\rm )} with hyperbolic vertex and
edge groups. Further{\rm ,} assume the two inclusions of $H$ are
quasi-isometric embeddings. Then $H$ is of finite height in $G$ if and
only if it is quasiconvex in $G$.
\end{theor}

The main theorem of \cite{GMRS} states:

\setcounter{section}{1}
\setcounter{theo}{0}
\begin{theor}[\!]
If $H$ is a quasiconvex subgroup of a hyperbolic group $G${\rm ,} then
$H$ has finite height.\label{gmrs}
\end{theor}

Thus the purpose of this paper is to prove the converse direction.

Certain group theoretic analogs of Thurston's combination theorems
\cite{ctm} were deduced in \cite{BF}. Extending the analogy with
\cite{ctm}, in this paper we prove quasiconvexity of certain surface
subgroups.

\setcounter{section}{5}
\setcounter{theo}{0}
\begin{propo}$\left.\right.$\vspace{.5pc}

\noindent Let $G = {G_1}*_H{G_2}$ be a hyperbolic group such that $G_1,
G_2, H$ are hyperbolic and the two inclusions of $H$ are
quasi-isometric embeddings. If $H$ is malnormal in one of $G_1$ or $G_2$
then $H$ is quasiconvex in $G$.
\end{propo}

The following corollary is a group-theoretic analog of a theorem
of Thurston's \cite{ctm}.

\setcounter{theo}{2}
\begin{coro}$\left.\right.$\vspace{.5pc}

\noindent Let $M_1$ be a hyperbolic atoroidal acylindrical $3$-manifold
and $S_1$ an incompressible surface in its boundary. Let $M_2$ be a
hyperbolic atoroidal $3$-manifold and $S_2$ an incompressible surface in
its boundary. If $S_1$ and $S_2$ are homeomorphic then gluing $M_1$ and
$M_2$ along this common boundary $S\, (= S_1 = S_2)$ one obtains a 
$3$-manifold $M$ such that
\begin{enumerate}
\renewcommand\labelenumi{\rm \arabic{enumi}.}
\item ${\pi_1}(M)$ is hyperbolic.

\item ${\pi_1}(S)$ is quasiconvex in ${\pi_1}(M)$.
\end{enumerate}
\end{coro}

\setcounter{section}{1}
\section{Preliminaries}

We start off with some preliminaries about hyperbolic metric spaces in
the sense of Gromov \cite{Gromov}. For details, see \cite{CDP,GhH}. Let
$(X,d)$ be a hyperbolic metric space.

\begin{definit}$\left.\right.$\vspace{.5pc}

\noindent A subset $Z$ of $X$ is said to be {\it $k$-quasiconvex} if any
geodesic joining $a,b\in Z$ lies in a $k$-neighborhood of $Z$. A subset
$Z$ is {\it quasiconvex} if it is $k$-quasiconvex for some $k$. A map
$f$ from one metric space $(Y,{d_Y})$ into another metric space
$(Z,{d_Z})$ is said to be a {\it $(K,\epsilon)$-quasi-isometric
embedding} if
\begin{equation*}
{\frac{1}{K}}({d_Y}({y_1},{y_2}))-\epsilon\leq{d_Z}(f({y_1}),
f({y_2}))\leq{K}{d_Y}({y_1},{y_2})+\epsilon.
\end{equation*}
If $f$ is a quasi-isometric embedding, and every point of $Z$ lies at a
uniformly bounded distance from some $f(y)$ then $f$ is said to be a
{\it quasi-isometry}. A $(K,{\epsilon})$-quasi-isometric embedding that
is a quasi-isometry will be called a $(K,{\epsilon})$-quasi-isometry.

A {\it $(K,\epsilon)$-quasigeodesic} is a $(K,\epsilon)$-quasi-isometric
embedding of a closed interval in~$\mathbb{R}$. A $(K,0)$-quasigeodesic
will also be called a $K$-quasigeodesic.
\end{definit}

\begin{definit}\cite{farb,Gromov2}$\left.\right.$\vspace{.5pc}

\noindent If $i: \Gamma_H \rightarrow \Gamma_G$ be an embedding of the
Cayley graph of $H$ into that of $G$, then the distortion function
is given by
\begin{equation*}
{\rm disto}(R) = {\rm Diam}_{\Gamma_H}({\Gamma_H}{\cap}B(R)),
\end{equation*}
where $B(R)$ is the ball of radius $R$ around $1\in{\Gamma_G}$.\vspace{.5pc}
\end{definit}

If $H$ is quasiconvex in $G$ the distortion function is linear and we
shall refer to $H$ as an undistorted subgroup. Else, $H$ will be termed
distorted. Note that the above definition makes sense for metric spaces
and their subspaces too.

\section{Trees of hyperbolic metric spaces}

For a general discussion of graphs of groups, see \cite{scott-wall}. In
this paper we will deal with graphs of hyperbolic groups satisfying the
quasi-isometrically embedded condition of \cite{BF}. We will need some
results from \cite{mitra3}.

\begin{definit}$\left.\right.$\vspace{.5pc}

\noindent A tree $(T)$ of hyperbolic metric spaces satisfying the q(uasi)
i(sometrically) embedded condition is a metric space $(X,d)$ admitting a
map $P : X \rightarrow T$ onto a simplicial tree $T$, such that there
exist $\delta{,} \epsilon$ and $K > 0$ satisfying the following:
\begin{enumerate}
\renewcommand\labelenumi{\arabic{enumi}.}
\item For all vertices $v\in{T}$, $X_v = P^{-1}(v) \subset X$ with the
induced path metric $d_v$ is a $\delta$-hyperbolic metric space.
Further, the inclusions ${i_v}:{X_v}\rightarrow{X}$ are uniformly
proper, i.e. for all $M > 0$, $v\in{T}$ and $x, y\in{X_v}$, there exists
$N > 0$ such that $d({i_v}(x),{i_v}(y)) \leq M$ implies ${d_v}(x,y) \leq
N$.

\item Let $e$ be an edge of $T$ with initial and final vertices $v_1$
and $v_2$ respectively. Let $X_e$ be the pre-image under $P$ of the
mid-point of $e$. Then $X_e$ with the induced path metric is
$\delta$-hyperbolic.

\item There exist maps ${f_e}:{X_e}{\times}[0,1]\rightarrow{X}$, such
that $f_e{|}_{{X_e}{\times}(0,1)}$ is an isometry onto the pre-image of
the interior of $e$ equipped with the path metric.

\item ${f_e}|_{{X_e}{\times}\{{0}\}}$ and
${f_e}|_{{X_e}{\times}\{{1}\}}$ are $(K,{\epsilon})$-quasi-isometric
embeddings into $X_{v_1}$ and $X_{v_2}$ respectively.
${f_e}|_{{X_e}{\times}\{{0}\}}$ and ${f_e}|_{{X_e}{\times}\{{1}\}}$ will
occasionally be referred to as $f_{v_1}$ and $f_{v_2}$ respectively.
\end{enumerate}
\end{definit}

$d_v$ and $d_e$ will denote path metrics on $X_v$ and $X_e$ respectively.
$i_v$, $i_e$ will denote inclusion of $X_v$, $X_e$ respectively into $X$.

We shall need a construction used in \cite{mitra3}. For convenience of
exposition, $T$ shall be assumed to be rooted, i.e. equipped with a base
vertex $v_0$. We shall refer to $X_{v_0}$ as $Y$. Let $v \neq v_0$ be a
vertex of $T$. Let $v_{-}$ be the penultimate vertex on the geodesic
edge path from $v_0$ to $v$. Let $e$ denote the directed edge from
${v_{-}}$ to $v$. Define ${\phi_v} :
{f_{e_{-}}}({X_{e_{-}}}{\times}\{{0}\}) \rightarrow
{f_{e_{-}}}({X_{e_{-}}}{\times}\{{1}\})$ as follows:

If $p{\in}{f_{e_{-}}}({X_e}{\times}\{0\}){\subset}{X_{v_{-}}}$, choose
$x\in{X_e}$ such that $p={f_{e_{-}}}(x{\times}\{0\})$ and define
\begin{equation*}
{\phi_v}(p) = {f_{e_{-}}}({x}{\times}\{{1}\}).
\end{equation*}

Note that in the above definition, $x$ is chosen from a set of bounded
diameter.

Let $\mu$ be a geodesic in $X_{v_{-}}$, joining $a, b \in
{f_{e_{-}}}({X_{e_{-}}}{\times}\{{0}\}) $. ${\Phi}_{v}({\mu})$ will
denote a geodesic in $X_v$ joining $\phi_v{(a)}$ and $\phi_v{(b)}$. Let
$X_{v_0} = Y$ and $ i = i_{v_0}$.

The next lemma follows easily from the fact that local quasigeodesics in
a hyperbolic metric space are quasigeodesics \cite{GhH}. If $x, y$ are
points in a hyperbolic metric space, $[x,y]$ will denote a geodesic
joining them.

\setcounter{theo}{0}
\begin{lem}
Given $\delta > 0${\rm ,} there exist $D, C_1$ such that if $a, b, c, d$
are vertices of a $\delta$-hyperbolic metric space $(Z,d)${\rm ,} with
${d}(a,[b,c])={d}(a,b)${\rm ,} ${d}(d,[b,c])={d}(c,d)$ and
${d}(b,c)\geq{D}$ then $[a,b]\cup{[b,c]}\cup{[c,d]}$ lies in a
$C_1$-neighborhood of any geodesic joining $a, d$.\label{perps}
\end{lem}

Given a geodesic segment $\lambda\subset{Y}$, we now recall from
\cite{mitra3} the construction of a quasi-convex set
$B_{\lambda}\subset{X}$ containing $i({\lambda})$.

\subsection*{\it Construction of quasiconvex sets}

Choose $C_2\geq{0}$ such that for all $e\in{T}$,
${f_e}({X_e}{\times}\{{0}\})$ and ${f_e}({X_e}{\times}\{{1}\})$ are
$C_2$-quasiconvex in the appropriate vertex spaces. Let $C{=}C_1{+}C_2$,
where $C_1$ is as in Lemma \ref{perps}.

For $Z\subset{X_v}$, let ${N_C}(Z)$ denote the $C$-neighborhood of $Z$,
that is the set of points at distance less than or equal to $C$ from
$Z$.

\begin{step}
{\rm Let $\mu\subset{X_v}$ be a geodesic segment in $({X_v},{d_v})$. Then
$P({\mu}) = v$. For each edge $e$ incident on $v$, but not lying on the
geodesic (in $T$) from $v_0$ to $v$, choose $p_e$, $q_e$ $\in
{N_C}({\mu}){\cap}{f_v}({X_e})$ such that ${d_v}({p_e},{q_e})$ is
maximal. Let ${v_1},{\ldots},{v_n}$ be terminal vertices of edges
${e_i}$ for which ${d_v}({p_{e_i}},{q_{e_i}}) > D$, where $D$ is as in
Lemma \ref{perps} above. Observe that there are only finitely many
$v_i$'s as $\mu$ is finite. Define
\begin{equation*}
{B^1}({\mu}) = {i_v}({\mu}){\cup}{\bigcup}_{k=1{\ldots} n}{{\Phi}_{v_i}}
({{\mu_i}}),
\end{equation*}
where $\mu_i$ is a geodesic in ${X_v}$ joining ${p_{e_i}}, {q_{e_i}}$.

Note that $P({B^1}{({\mu})})\subset{T}$ is a finite tree.

The reason for insisting that the edges $e$ do not lie on the geodesic
from $v_0$ to $v$ is to prevent `backtracking' in Step 2 below.}
\end{step}

\begin{step}
{\rm Step 1 above constructs $B^1({\lambda})$ in particular. We proceed
inductively. Suppose that $B^m({\lambda})$ has been constructed such
that the convex hull of $P({B^m({\lambda})}) \subset {T}$ is a finite
tree. Let $ \{{w_1}, \ldots , {w_n}\} =
P({B^m({\lambda})}){\setminus}P({B^{m-1}}({\lambda}))$. (Note that $n$
may depend on $m$, but we avoid repeated indices for notational
convenience.) Assume further that
${P^{-1}}({v_k}){\cap}{B^m({\lambda})}$ is a path of the form
${i_{v_k}}({\lambda}_{k})$, where $\lambda_k$ is a geodesic in
$({X_{v_k}},{d_{v_k}})$. Define
\begin{equation*}
{B^{m+1}}({\lambda}) = {B^m}({\lambda}){\cup}{\bigcup}_{k=1\ldots{n}}
(B^1({\lambda_k})),
\end{equation*}
where $B^1({\lambda_k})$ is defined in Step 1 above.

Since each $\lambda_k$ is a finite geodesic segment in $\Gamma_H$, the
convex hull of $P({B^{m+1}}{\lambda})$ is a finite subtree of $T$.
Further, $P^{-1}{(v)}{\cap}B^{m+1}({\lambda})$ is of the form
${i_v}({\lambda_v})$ for all $v\in{P({B^{m+1}({\lambda})})}$. This
enables us to continue inductively. Define 
\begin{equation*}
B({\lambda}) = {\cup}_{m\geq{0}}B^m{\lambda}.
\end{equation*}

Note that the convex hull of $P({B({\lambda})})$ in $T$ is a locally
finite tree $T_1$. Further $B({\lambda})\cap{{P^{-1}}(v)}$ is a geodesic
in $X_v$ for $v\in{T_1}$ and is empty otherwise.}
\end{step}

\subsection*{\it Construction of retraction}

One of the main theorems of \cite{mitra3} states that $B({\lambda})$
constructed above is uniformly quasiconvex. To do this we constructed a
retraction $\Pi_{\lambda} $ from (the vertex set of) $X$ onto
${B_\lambda}$ and showed that there exists $C_0\geq{0}$ such that
${d_X}({\Pi_\lambda}(x),{\Pi_\lambda}(y))\leq C_0{d_X}(x,y)$. Recall
this construction from \cite{mitra3}. Let $\pi_v :
X_v\rightarrow{\lambda_v}$ be a nearest point projection of $X_v$ onto
$\lambda_v$. $\Pi_\lambda$ is defined on ${\bigcup}_{v\in{T_1}}X_v$ by
\begin{equation*}
\Pi_{\lambda}(x) = {i_v}{\cdot}{\pi_v}(x)\ \ \hbox{for}\ \ x\in{X_v}.
\end{equation*}

If $x\in{P^{-1}}(T\setminus{T_1})$ choose $x_1\in{P^{-1}}({T_1})$
such that ${d}(x,{x_1}) = {d}(x,{P^{-1}}({T_1}))$ and define
${\Pi_{\lambda}^{\prime}}(x) = x_1$. Next define
$\Pi_\lambda{(x)} = {\Pi_\lambda}\cdot{\Pi_{\lambda}^{\prime}(x)}$.

\begin{theor}[\!]\hskip -.5pc {\bf \cite{mitra3}.}\ \ 
There exists $C_0\geq{0}$ such that
${d}({\Pi_\lambda}(x),{\Pi_\lambda}(y))\leq C_0{d}(x,y)$ for $x, y$
vertices of $X$. Further{\rm ,} $B({\lambda})$ is $C_0$-quasiconvex.
\label{mainref}
\end{theor}

We need one final lemma from \cite{mitra3}. Let $i: Y \rightarrow X$
denote inclusion.

\begin{lem}
There exists $A > 0${\rm ,} such that if
$a\in{P^{-1}}(v){\cap}{B{({\lambda})}}$ for some $v\in{T_1}$ then there
exists $b\in{i({\lambda})} = P^{-1}{(v_0)}\cap{B{({\lambda})}}$ with
${d}(a,b) \leq {A}{d_T}(Pa,Pb)$. Further{\rm ,} let $v_0 , v_1 , \ldots,
v_n =v$ be the sequence of vertices on a geodesic in $T$ connecting the
root vertex $v_0$ to $v$. There exists a sequence $b = a_0 , a_1 ,
\ldots, a_n = a$ with $a_i \in P^{-1}{(v_i)}\cap{B({\lambda}})$ such
that $d({a_i},{a_j}) \leq A{d_T} (P{a_i},P{a_j}) = A{d_T}({v_i},{v_j})$.
\label{connectionlemma}
\end{lem}

The above lemma says that we can construct a quasi-isometric section of
a geodesic segment $[{v_0},v]$ ending at $a$.

\begin{definit}$\left.\right.$\vspace{.5pc}

\noindent An {\it $A$-quasi-isometric section} of $[{v_0},v]$
ending at $a \in P^{-1}(v){\cap}{B({\lambda})}$ is a sequence of points
in $X$ satisfying the conclusions of Lemma \ref{connectionlemma} above.\vspace{.5pc}
\end{definit}

Note that the quasi-isometric sections considered are all images of
$[{v_0},v]$ where $v_0$ is the root vertex of $T$. Abusing notation
slightly we will refer to the map or its image as a quasi-isometric
section.

So far we have considered a tree of hyperbolic metric spaces.
It is time to introduce the relevant groups.

Let $G$ be a hyperbolic group acting cocompactly on a simplicial tree
$T$ such that all vertex and edge stabilizers are hyperbolic. Also
suppose that every inclusion of an edge stabilizer in a vertex
stabilizer is a quasi-isometric embedding. Let $\mathcal G$ denote the
quotient graph $T/G$. The metric on $T$ will be denoted by $d_T$. Assume
$\mathcal{G}$ has only one edge and $H$ is the stabilizer of this edge.
This is the situation when $G$ splits over $H$.

Suppose $H$ is a vertex or edge subgroup. Further, suppose $H$ is
distorted in $G$. We would like to show that $H$ has infinite height.
Here is a brief sketch of the proof of the main theorem of this paper:

Since $H$ is distorted, there exist geodesics $\lambda_i \subset
\Gamma_H$ such that geodesics in $\Gamma_G$ joining the end points of
$\lambda_i$ leave larger and larger neighborhoods of $\Gamma_H$. From
the construction of $B({\lambda})$ it follows that the diameters ${\rm
dia}({P(B({\lambda_i}))})\rightarrow\infty$ as $i\rightarrow\infty$. The
edges of $T$ can be lifted to $\Gamma_G$ and one can after a pigeon-hole
principle argument look upon these lifts as conjugating elements. The
geodesics in $B({\lambda_i})\cap{P^{-1}}(v)$ can be thought of as
elements of $H$. Thus as $i\rightarrow\infty$ one obtains a sequence of
elements $g_i\in{G}$ such that $\cap{g_i^{-1}}H{g_i} \neq 1$. This
proves that $H$ has infinite height. The next section is devoted to
making this rigorous.

\section{Proof of Main Theorem}

We start our discussion with a basic lemma.

\setcounter{theo}{0}
\begin{lem}
If $X_{v_0} = Y$ is distorted in $X${\rm ,} there exist a sequence of
geodesics $\lambda_i$ in $Y$ such that ${\rm dia}(P(B({\lambda_i})))
\rightarrow \infty$ as $i \rightarrow \infty${\rm ,} where the diameter
is calculated with respect to the metric $d_T$.\label{dia}
\end{lem}

\begin{proof}
It follows from Lemma \ref{connectionlemma} that $B({\lambda_i})$ lies
in an $A$ dia$(P(B({\lambda_i})))$ neighborhood of $i({\lambda_i})$ and
hence of $Y$. Further from Theorem \ref{mainref} a geodesic in $X$
joining the end points of $i({\lambda_i})$ lies in a (uniform)
$C_0$-neighborhood of $B({\lambda_i})$.

Since $Y$ is distorted in $X$, there exist $\lambda_i \subset Y$ such
that geodesics in $X$ joining end points of $\lambda_i$ leave an 
$i$-neighborhood of $Y$ for $i = 1, 2, \ldots$.

Hence $i \leq A$ dia $(P(B({\lambda_i}))) +C$.

The lemma follows.\hfill \ab
\end{proof}

\subsection*{\it Construction of hallways}

We would like to construct certain special subsets of $B({\lambda})$
closely related to the essential hallways of Bestvina and Feighn
\cite{BF}. We retain the terminology.

\begin{definit}$\left.\right.$\vspace{.5pc}

\noindent A disk $f : [0,m]{\times}{I} \rightarrow X$ is a hallway of
length $m$ if it satisfies:

\begin{enumerate}
\renewcommand\labelenumi{\arabic{enumi}.}
\item $f^{-1} ({\cup}{X_v} : v \in T) = \{0, 1, \ldots, m \}{\times} I$.

\item $f$ maps $i{\times}I$ to a geodesic in $X_v$ for some vertex
space.

\item $( P\circ{f} ) : [0,m]{\times}{I} \rightarrow T$ factors through
the canonical retraction to $[0,m]$ and an isometry of $[0,m]$ to $T$.
\end{enumerate}
\end{definit}

\begin{definit}$\left.\right.$\vspace{.5pc}

\noindent A hallway is $\rho$-thin if $d({f(i,t)},{f({i+1},t)}) \leq
\rho$ for all $i, t$.\vspace{.5pc}
\end{definit}

We will now construct $A$-thin hallways using the quasi-isometric
sections of Lemma \ref{connectionlemma}. The arguments are carried out
for trees of metric spaces.

Given $\lambda$ and $x\in{B({\lambda})}$ let $\Sigma^x_{\lambda}$ be an
$A$-quasi-isometric section of $[{v_0},P(x)]$ into $B({\lambda})$ ending
at $x$. From Lemma \ref{connectionlemma} such quasi-isometric sections
exist. Further, if $a \in \Sigma^x_{\lambda}$ then define
$\sigma^x_{\lambda} (a)$ to be a point
$i({\lambda}){\cap}\Sigma^x_{\lambda}$. The choice involved in the
definition of $\sigma^x_{\lambda} (a)$ is bounded purely in terms of
$A$. 

\begin{lem}
Suppose $Y = X_{v_0}$ is distorted in $X$. Then there exist geodesics
$\lambda_i\subset{Y}, {a_i}, {b_i}, {x_i}, {y_i} \in B({\lambda_i})$
such that
\begin{enumerate}
\renewcommand\labelenumi{\rm \arabic{enumi}.}
\item $d({x_i},{y_i})\leq 1$.

\item $P({x_i}) = P({y_i})$.

\item $\mu_i$ is a geodesic subsegment of $\lambda_i$ in $Y$ joining
$\sigma^{a_i}_{\lambda_i} ({x_i})$ and $\sigma^{b_i}_{\lambda_i}
({y_i})$ with length of $\mu_i$ greater than or equal to $i$.
\end{enumerate}\label{hallways}
\end{lem}

\begin{proof}
Suppose not. Then there exists $C \geq 0$ such that for all geodesics
$\lambda_i$ in $Y$ and all ${a_i}, {b_i}, {x_i}, {y_i} \in
B({\lambda_i})$ satisfying

\begin{enumerate}
\renewcommand\labelenumi{\rm \arabic{enumi}.}
\item ${a_i}, {b_i}, {x_i}, {y_i} \in B({\lambda_i})$.

\item $d({x_i},{y_i})\leq 1$.

\item $P({x_i}) = P({y_i})$.

\item $\mu_i$ is a geodesic subsegment of $\lambda_i$ in $Y$ joining
$\sigma^{a_i}_{\lambda_i} ({x_i})$ and $\sigma^{b_i}_{\lambda_i}
({y_i})$.
\end{enumerate}

We have length of $\mu_i$ less than or equal to $C$. For all
$x\in{B({\lambda_i})}$ choose $a \in B({\lambda_i})$ such that $x \in
\Sigma^a_{\lambda_i}$ and define
\begin{equation*}
\pi (x) = \sigma^a_{\lambda_i}{(x)}.
\end{equation*}

Recall that $\pi{(x)}$ is chosen from a set of (uniformly) bounded
diameter. Thus we might as well take $ a = x$. Note that $\pi$ defines a
retraction of $B({\lambda_i})$ onto $\lambda_i$.

For any $x, y\in B({\lambda_i})$ such that $P(x) = P(y)$ we have
$d({\pi}(x),{\pi}(y)) \leq C d(x,y)$.

Next suppose $x, y \in B({\lambda_i})$, $d(P(x),P(y)) = 1$ and $d(x,y)
\leq A$. Assume without loss of generality $d(P(x),{v_0}) <
d(P(y),{v_0})$. Then by Lemma \ref{connectionlemma} there exists
$z\in{B({\lambda_i})}{\cap}{\Sigma^y_{\lambda_i}}$ such that $P(x) =
P(z)$, $d(x,z) \leq 2A$ and hence $d({\pi}(x),{\pi}(y)) \leq 2AC + C$.

Hence there exists $C^{\prime}$ such that for any $\lambda_i$ and $x, y
\in \lambda_i$, $d({\pi}(x),{\pi}(y)) \leq C^{\prime}d(x,y)$. Thus
$\lambda_i$ is uniformly quasiconvex in $B({\lambda_i})$ and hence (by
Theorem \ref{mainref}) in $X$.

Therefore $Y$ is quasiconvex in $X$, contradicting the hypothesis.
\hfill\ab
\end{proof}

\begin{definit}$\left.\right.$\vspace{.5pc}

\noindent An $A$-thin hallway $\mathcal{H}$ with ends $\mu_0, \mu_n$
trapped by $A$-quasi-isometric sections $\Sigma_1$ and $\Sigma_2$ is a
collection of geodesics $\mu_i \subset X_{v_i}, i = 0, \ldots, n$ such
that

\begin{enumerate}
\renewcommand\labelenumi{\rm \arabic{enumi}.}
\item $v_0, {\ldots}, v_n$ are successive vertices on a geodesic
$[{v_0}, {v_n}]$ in $T$.

\item $\mu_i$ joins $\Sigma_1 ({v_i})$ to $\Sigma_2 ({v_i})$.
\end{enumerate}

As before $n$ is called the length of the hallway.
\end{definit}

Note that the geodesics are allowed to have length 0.

\begin{coro}{\rm Existence of hallways} $\left.\right.$\vspace{.5pc}

\noindent Suppose $Y$ is distorted in $X$.
Then there exist geodesics $\lambda_i \subset Y$ and $A$-thin hallways
${\mathcal{H}}_i$ with ends $\lambda_i, \eta_i$ trapped by
quasi-isometric sections $\Sigma_{1i}, \Sigma_{2i}$ such that the
lengths of $\lambda_i$ and the hallway ${\mathcal{H}}_i$ are greater
than $i$.\label{hallways2}
\end{coro}

\begin{proof}
From Lemma \ref{hallways} there exist geodesics $\lambda_i\subset{Y}$,
${a_i}, {b_i}, {x_i}, {y_i} \in B({\lambda_i})$ such that

\begin{enumerate}
\renewcommand\labelenumi{\rm \arabic{enumi}.}
\item $d({x_i},{y_i})\leq 1$.

\item $P({x_i}) = P({y_i})$.

\item $\mu_i$ is a geodesic subsegment of $\lambda_i$ in $Y$ joining
$\sigma^{a_i}_{\lambda_i} ({x_i})$ and $\sigma^{b_i}_{\lambda_i}
({y_i})$ with length of $\mu_i$ greater than $i$.
\end{enumerate}

Take $\Sigma_{1i} = \Sigma^{a_i}_{\lambda_i}$ $\Sigma_{2i} =
\Sigma^{b_i}_{\lambda_i}$ and rename $\mu_i$ as $\lambda_i$ (we are
abusing notation slightly here).

Passing to a subsequence if necessary and arguing as in Lemma \ref{dia}
we can assume that the length of ${\mathcal{H}}_i$ is greater than $i$.

The corollary follows.\hfill \ab
\end{proof}

\subsection*{\it Construction of annuli}

The discussion so far has not entailed the use of group actions. We
would like to establish a dictionary between the geometric objects
constructed above and elements of a group $G$ acting on $T$.

Let $G$ be a hyperbolic group acting cocompactly on a simplicial tree
$T$ such that all vertex and edge stabilizers are hyperbolic. Also
suppose that every inclusion of an edge stabilizer in a vertex
stabilizer is a quasi-isometric embedding. Let $\mathcal G$ denote the
quotient graph $T/G$. The metric on $T$ will be denoted by $d_T$. Assume
$\mathcal{G}$ has only one edge and $H$ is the stabilizer of this edge.
This is the situation when $G$ splits over $H$. That is $G =
{G_1}*_H{G_2}$ or $G = {G_1}*_H$. Then by the restrictions on the
$G$-action on $T$, the inclusions of $H$ into $G_i$ are quasi-isometric
embeddings.

The stabilizers of edges of $T$ are conjugates of $H$. We can take
$\Gamma_H = {X_{v_0}} = Y$, $\Gamma_G = X$ and $i : Y \rightarrow X$ the
natural inclusion. Let $\lambda \subset Y$ be a geodesic.

Recall the construction of $B({\lambda})$ from the previous section.
$B({\lambda})$ was constructed as the union of certain geodesics
$\lambda_i \subset X_{v_i}$. Further, each $\lambda_i$ was in the image
of an edge space. Therefore if $\lambda_i$ has $a_i$, $b_i$ as its
end points, then $a_i^{-1}{b_i} \in H$.

We need to now examine the hallways constructed above. Let
${\mathcal{H}} = {\cup}_{i=0, \ldots, n}\mu_i$ be an $A$-thin hallway
trapped between quasi-isometric sections $\Sigma_1$ and $\Sigma_2$ with
ends $\mu_0$ and $\mu_n$. Note that each $\mu_i$ is a geodesic
subsegment of some $\lambda_i$ joining $a_i$, $b_i$ and $a_i^{-1}{b_i}
\in H$.

Since edge spaces are (uniformly) quasi-isometrically embedded in vertex
spaces, there exists a constant $D_1$ such that if $\mu_i$ joins $c_i$,
$d_i$ then ${c_i}^{-1}{d_i} = {u_i}{h_i}{v_i}$, where $h_i \in H$,
${|}u_i{|} \leq {\frac{D_1}{2}}$ and ${|}v_i{|} \leq {\frac{D_1}{2}}$.
($|.|$ denotes length.) Also, from the definition of $A$-thin hallways
trapped between quasi-isometric sections, we have
\begin{align*}
&{\Sigma_1}(i) = c_i,\\[.2pc]
&{\Sigma_2}(i) = d_i,\\[.2pc]
&|{\Sigma_1}(i)^{-1}{\Sigma_2}(i)| \leq {\frac{D_1}{2}} \ \ {\rm for\ 
all} \ \ i.
\end{align*}

\begin{definit}$\left.\right.$\vspace{.5pc}

\noindent An $(A+{D_1})$-thin $H$-hallway $\mathcal{H}$ with ends
$\mu_0$, $\mu_n$ trapped by $(A+{D_1})$-quasi-isometric sections
$\Sigma_1$ and $\Sigma_2$ is a collection of geodesics $\mu_i \subset
X_{v_i}$, $i = 0, \ldots, n$ such that

\begin{enumerate}
\renewcommand\labelenumi{\arabic{enumi}.}
\item $v_0 , {\ldots}, v_n$ are successive vertices on a geodesic
$[{v_0},{v_n}]$ in $T$.

\item $\mu_i$ joins $\Sigma_1 ({v_i}) = c_i$ to $\Sigma_2 ({v_i}) =
d_i$.

\item $c_i^{-1}{d_i} \in H$.
\end{enumerate}
\end{definit}

The following lemma is the group-theoretic counterpart of Corollary
\ref{hallways2} and follows from the discussion above.

\begin{lem}
Suppose $Y (= X_{v_0} = \Gamma_H)$ is distorted in $X (= \Gamma_G )$.
Then there exist geodesics $\lambda_i \subset Y$ and $(A+{D_1})$-thin
$H$-hallways ${\mathcal{H}}_i$ with ends $\lambda_i${\rm ,} $\eta_i$
trapped by $A+{D_1}$-quasi-isometric sections $\Sigma_{1i}${\rm ,}
$\Sigma_{2i}$ such that the lengths of $\lambda_i$ and the hallway
${\mathcal{H}}_i$ are greater than~$i$.\label{H-hallways}
\end{lem}

We would now like to paste two of these $H$-hallways together along a
common bounding quasi-isometric section.

Given $n > 0$ consider $(A+{D_1})$-thin hallways ${\mathcal{H}}_i$ with
one end $\lambda_i \subset Y = \Gamma_H$ of length $n$. Clearly there
exist infinitely many distinct such from Lemma \ref{H-hallways} (taking
a long enough hallway with one end in $Y$ and truncating it to one of
length $n$ gives such a hallway).

\begin{definit}$\left.\right.$\vspace{.5pc}

\noindent The ordered boundary $\Delta_{\mathcal{H}}$ of an $H$-hallway
$\mathcal{H}$ of length $n$ trapped by quasi-isometric sections
$\Sigma_1$, $\Sigma_2$ is given by
\begin{equation*}
\Delta_{\mathcal{H}} = \{ {\Sigma_1}({v_{j-1}})^{-1}{\Sigma_1}({v_j}),
{\Sigma_2}({v_{j-1}})^{-1}{\Sigma_1}({v_j}), : j = 1 \ldots n \},
\end{equation*}
where $[v_0,v_n] \subset T$ is the geodesic in $T$ to which
$\mathcal{H}$ maps under $P$.

The $i$th element of the above set will be denoted by
$\Delta_{\mathcal{H}} (i)$.

If the hallway is $A+{D_1}$-thin, then
$|{\Sigma_i}({v_{j-1}})^{-1}{\Sigma_i}({v_j})| \leq A+D_1$.

Since there exist infinitely many distinct $(A+{D_1})$-thin $H$-hallways
of length $n$ and only finitely many words in $G$ of length less than or
equal to $(A+{D_1})$, there exist (by the pigeon-hole principle)
infinitely many distinct $H$-hallways of length $n$ with the same
ordered boundary $\Delta$.

Choose two such hallways and glue one to the `reflection' of the other.
More precisely, let ${\mathcal{H}}_i = \cup_{j=1\ldots n} \mu_{ij}$ for $i =
1, 2$ be two such hallways. Let $\mu_{ij}$ have $a_{ij}, b_{ij} \in
X_{v_j} \subset \Gamma_G$ as its end points.

Then since ${\mathcal{H}}_i$ are $(A+{D_1})$-thin $H$-hallways with the
same ordered boundary, we have
\begin{align*}
&a_{ij}^{-1}b_{ij} \in H,\\[.2pc]
&a_{1j}^{-1}a_{1,j+1} = a_{2j}^{-1}a_{2,j+1},\\[.2pc]
&b_{1j}^{-1}b_{1,j+1} = b_{2j}^{-1}b_{2,j+1}.
\end{align*}

Let $\eta_j$ denote a geodesic in $X_{v_j}$ joining $a_{1j}$ and $c_{1j}
= b_{1j}b_{2j}^{-1}{a_{2j}}$. Then $\mathcal{H} = \cup_{j=1\ldots n}
\eta_j$ is an $(A+{D_1})$-thin $H$-hallway. If $\Delta$ be its ordered
boundary, then it follows from the above equations that $\Delta{(2i)} =
{\Delta}(2i-1)$ for $i = 1{\ldots}n$.
\end{definit}

\begin{definit}$\left.\right.$\vspace{.5pc}

\noindent An $H$-hallway of length $n$ with ordered boundary $\Delta$ is
called an $H$-annulus if $\Delta{(2i)}={\Delta}(2i-1)$ for $i =
1{\ldots}n$.\vspace{.5pc}
\end{definit}

The above definition is related to the annuli of Bestvina and Feighn
\cite{BF}.

From the above discussion and Lemma \ref{H-hallways} the following
crucial theorem follows:

\begin{theor}[\!]
Suppose \hbox{$Y (=X_{v_0} = \Gamma_H)$} is distorted in \hbox{$X (= \Gamma_G)$}.
Then there exist geodesics $\lambda_i \subset Y$ and $(A+{D_1})$-thin
$H$-annuli ${\mathcal{H}}_i$ with ends $\lambda_i${\rm ,} $\eta_i$
trapped by $(A+{D_1})$-quasi-isometric sections $\Sigma_{1i}${\rm ,}
$\Sigma_{2i}$ such that the lengths of $\lambda_i$ and the hallway
${\mathcal{H}}_i$ are greater than $i$. In fact there exist infinitely
many distinct such $H$-annuli with the same ordered boundary.
\label{H-annuli}
\end{theor}

The main theorem of this paper follows from Theorem \ref{H-annuli} by
unravelling definitions. We state this below.

\begin{theor}[\!]
Let $G$ be a hyperbolic group splitting over $H$ {\rm (}i.e. $G =
{G_1}{*_H}{G_2}$ or $G = {G_1}{*_H}${\rm )} with hyperbolic vertex and
edge groups. Further{\rm ,} assume the two inclusions of $H$ are
quasi-isometric embeddings. Then $H$ is of finite height in $G$ if and
only if it is quasiconvex in $G$.\label{main}
\end{theor}

\begin{proof}
Suppose $H$ is distorted in $G$. Then from Theorem \ref{H-annuli} there
exists an $H$-annulus $\mathcal{H} = \cup_{i=0\ldots n}\lambda_i$ of
length $n$ such that $|{\lambda_0}| > n$. (In fact there are infinitely
many distinct such. However, we start off with one in the interests of
notation.)

Let $\Delta$ be the ordered boundary of $\mathcal{H}$. By definition of
$H$-annulus $\Delta{(2i)} = {\Delta}(2i-1)$ for $i = 1{\ldots}n$. Let
$c_i, d_i$ be the endpoints of $\lambda_i$ such that
\begin{equation*}
{\Delta}(2i-1) = c_{i-1}^{-1}c_i = d_{i-1}^{-1}d_i = {\Delta}(2i).
\end{equation*}

Also $c_i^{-1}d_i = h_i \in H$. Let $g_i = {\Delta}(2){\ldots}
{\Delta}(2i)$. Reading relations around `quadrilaterals' we have,
\begin{equation*}
h_{i-1} = {\Delta}(2i){h_i}{\Delta}(2i)^{-1}\ \ \hbox{for all}\ \ i = 1
\ldots n.
\end{equation*}

Therefore
\begin{equation*}
h_{0} = {g_i}{h_i}{g_i}^{-1}\ \ \hbox{for all}\ \ i = 1 \ldots n.
\end{equation*}

Recall that $P : \Gamma_G \rightarrow T$ is the projection onto $T$.
Since $P({c_0}{g_i}) \neq P({c_0}{g_j})$ for $i \neq j$ we have $n$
essentially distinct conjugates ${g_i}H{g_i}^{-1}$ whose intersection
contains $h_0 \neq 1$.

Now we need the fact that there are infinitely many distinct $H$-annuli
(Theorem \ref{H-annuli}) with the same ordered boundary. Without loss of
generality, let this boundary be $\Delta$ above. The above argument then
furnishes infinitely many distinct $h \in
{H}\cap_{i=1\ldots n}{g_i}H{g_i}^{-1}$.

Thus given any $n > 0$ there exist $n + 1$ essentially distinct
conjugates of $H$ whose intersection is infinite. Therefore $H$ has
infinite height. Along with Theorem~\ref{gmrs} this proves the
Theorem.\hfill \ab
\end{proof}

\section{Consequences and questions}

\subsection*{\it Malnormality}

We deduce a couple of group-theoretic consequences of Theorem~\ref{main}.

\begin{definit}$\left.\right.$\vspace{.5pc}

\noindent A subgroup $H$ of a group $G$ is said to be {\it malnormal} in
$G$ if $gH{g^{-1}}{\cap}H = 1$ for all $g \notin H$.
\end{definit}

\setcounter{theo}{0}
\begin{propo}$\left.\right.$\vspace{.5pc}

\noindent Let $G = {G_1}*_H{G_2}$ be a hyperbolic group such that $G_1,
G_2, H$ are hyperbolic and the two inclusions of $H$ are quasi-isometric
embeddings. If $H$ is malnormal in one of $G_1$ or $G_2$ then $H$ is
quasiconvex in $G$. \label{maln}
\end{propo}

\begin{proof}
Assume without loss of generality that $H$ is malnormal in $G_2$. Let $g
\in G{\setminus}H$ and $h, {h_1}\in{H}$ be such that $ghg^{-1} = h_1
\neq 1$. Let $g = {a_1}{b_1}{\ldots}{a_n}{b_n}$ with $a_i \in G_1$ and
$b_i \in G_2$. Then by normal form for free products with amalgamation
(\cite{ls}, p.~178) we have ${b_n}H{b_n^{-1}} \in H$ and hence $b_n \in
H$ by malnormality of $H$ in $G_2$. Continuing inductively, we get
${a_i}{\ldots}{a_n}h{a_n^{-1}}{\ldots}{a_i^{-1}}$ and $b_i \in H$ for
all $i = 1{\ldots}n$. In particular $g{\in}G_1$. Therefore
$H{\cap}{gHg^{-1}} \neq 1$ implies $g{\in}G_1$.

Since $H$ is quasi-isometrically embedded in $G_1$ we have by Theorem
\ref{gmrs} that $H$ has finite height in $G_1$. Therefore by the above
argument $H$ has finite height in $G$. Finally by Theorem~\ref{main}, $H$
is quasiconvex in $G$.\hfill \ab
\end{proof}

The above proposition holds good if malnormal is replaced by height
zero.

A similar argument using Britton's lemma~(\cite{ls}, p.~178) gives the
following:

\begin{propo}$\left.\right.$\vspace{.5pc}

\noindent Let $G = {G_1}*_H$ be a hyperbolic group such that $G_1, H$
are hyperbolic and the two images $H_1, H_2$ of $H$ are quasiconvex in
$G_1$. If $g{H_1}g^{-1}{\cap}{H_2}$ is finite for all $g \in G_1$ then
$H$ is quasiconvex in $G$. \label{malhn}
\end{propo}

The hypotheses in the above propositions cannot be relaxed as the
following example shows.

\begin{examm}
Let $G_i = \{ {a_i},{b_{1i}},{b_{2i}},{c_{1i}},{c_{2i}} |
{a_i}{b_{ji}}{a_i^{-1}} = c_{ji} , j = 1, 2 \}$ be two copies (for $i =
1, 2$) of a group isomorphic to the free group on 3 generators.

Let $H = \{ {b_1}, b_2 , c_1 , c_2 \}$ be the free group on 4
generators. Let $i : H \rightarrow G_1$ be given by sending $b_i$ to
$b_{i1}$ and $c_i$ to $c_{i1}$ for $i = 1, 2$.

Let $j : H \rightarrow G_2$ be given by sending $b_i$ to $b_{i2}$ for $i
= 1, 2$ and $c_i$ to `long words' $u_i$ in $c_{12}$ and $c_{22}$ such
that the `flare' condition of \cite{BF} is satisfied for the free
product with amalgamation $G = {G_1}{*_H}{G_2}$.

In fact one gets 
\begin{equation*}
G = \langle {a_1}, {a_2}, {c_{1}}, {c_{2}} |
{{a_1}{a_2^{-1}}{c_i}{a_2}{a_1^{-1}} = {u_i}({c_1},{c_2})}, i = 1, 2
\rangle
\end{equation*}
such that this is a small cancellation presentation with $G$ hyperbolic.

It is clear that the subgroup generated by ${c_1}, {c_2}$ is a free
group on two generators with infinite height in $G$. Hence the
amalgamating subgroup $H$ above is of infinite height.
\end{examm}

In \cite{ctm} McMullen shows that glueing an acylindrical, atoroidal
hyperbolic 3-manifold to another hyperbolic atoroidal 3-manifold 
along a common incompressible boundary surface $S$ gives a hyperbolic
3-manifold in which $S$ is quasifuchsian. We deduce the following
group theoretic version of this from Proposition~\ref{maln} above.

\begin{coro}$\left.\right.$\vspace{.5pc}

\noindent Let $M_1$ be a hyperbolic atoroidal acylindrical $3$-manifold
and $S_1$ an incompressible surface in its boundary. Let $M_2$ be a
hyperbolic atoroidal $3$-manifold and $S_2$ an incompressible surface in
its boundary. If $S_1$ and $S_2$ are homeomorphic then glueing $M_1$ and
$M_2$ along this common boundary $S\, (= S_1 = S_2)$ one obtains a
$3$-manifold $M$ such that

\begin{enumerate}
\renewcommand\labelenumi{\rm \arabic{enumi}.}
\item ${\pi_1}(M)$ is hyperbolic.

\item ${\pi_1}(S)$ is quasiconvex in ${\pi_1}(M)$.
\end{enumerate}\label{ctm}\vspace{-1pc}
\end{coro}

\begin{proof}
Hyperbolicity of $\pi_1 (M)$ follows from the combination theorem of
Bestvina and Feighn \cite{BF}. Quasiconvexity follows from Proposition
\ref{maln} above.\hfill \ab
\end{proof}

Using Proposition \ref{malhn} one can deduce similar results.

\subsection*{\it Graphs of hyperbolic groups}

The main argument of this paper does not generalize directly to graphs
of hyperbolic groups satisfying the quasi-isometrically embedded
condition. Given a distorted edge or vertex group $H \subset G$, the
pigeon-hole principle argument of the previous section does furnish an
edge group $H_1$ of infinite height in $G$ such that a conjugate of $H$
intersects $H_1$ in a distorted subgroup of $G$.

However $H$ and $H_1$ need not be the same. The basic problem lies in
dealing with quasiconvex subgroups of edge (or vertex) groups that are
distorted in $G$. We state the problem explicitly:

\begin{quest}
Suppose $G$ splits over $H$ satisfying the hypothesis of 
Theorem~\ref{main} and $H_1$ is a quasiconvex subgroup of $H$. If $H_1$ has
finite height in $G$ is it quasiconvex in $G$? More generally, if $H_1$
is an edge group in a hyperbolic graph of hyperbolic groups satisfying
the qi-embedded condition, is $H$ quasiconvex in $G$ if and only if it
has finite height in $G$?\vspace{.5pc}
\end{quest}

The above question is a special case of the general question of Swarup
on characterizing quasiconvexity in terms of finiteness of height.

There are two cases where a complete answer to the above question is
known. These are extensions of $\mathbb{Z}$ by surface groups
\cite{scottswar} or free groups \cite{BFH,mitra4}. Both these solutions
involve a detailed analysis of the ending laminations \cite{mitra2}.

\subsection*{\it Other questions}

A closely related problem \cite{bestvinahp,mitra} can be
formulated in more geometric terms:

\begin{quest}
Let $X_G$ be a finite 2 complex with fundamental group $G$. Let $X_H$ be
a cover of $X_G$ corresponding to the finitely presented subgroup $H$.
Let $I(x)$ be the injectivity radius of $X_H$ at $x$.

Does $I(x) \rightarrow \infty$ as $x \rightarrow \infty$ imply that
$H$ is quasi-isometrically embedded in $G$?
\end{quest}

A positive answer to this question for $G$ hyperbolic would provide
a positive answer to Swarup's question.

The answer to this question is negative if one allows $G$
to be only finitely generated instead of finitely presented
as the following example shows:

\begin{examm}
Let $F = \{a, b, c, d \}$ denote the free group on four generators. Let
${u_i} = ab^i $ and ${v_i} = cd^{f(i)}$ for some function $f :
\mathbb{N} \rightarrow \mathbb{N}$. Introducing a stable letter $t$
conjugating $u_i$ to $v_i$ one has a finitely generated HNN extension
$G$. The free subgroup generated by $a, b$ provides a negative answer to
the question above for suitable choice of $f$. In fact one only requires
that $f$ grows faster than any linear function. 

If $f$ is recursive one can embed the resultant $G$ in a finitely
presented group by Higman's embedding theorem. But then one might lose
malnormality of the free subgroup generated by $a, b$. A closely related
example was shown to the author by Steve Gersten.
\end{examm}

A counterexample to the general question of Swarup might provide a means
of constructing acyclic non-hyperbolic finitely presented groups without
$(Z+Z)$ answering a question of Bestvina and Brady \cite{bestvinahp}.
Suppose $H$ is a malnormal torsion-free hyperbolic subgroup of a
hyperbolic torsion-free group $G$. If $H$ is distorted in $G$, then
doubling $G$ along $H$ (i.e. $G{*_H}G$) one gets a finitely presented
acyclic group which is not hyperbolic, nor does it contain $(Z+Z)$. This
was independently observed by Sageev.

On the other hand one might develop an analog of Thurston's theory of
pleated surfaces \cite{Thurstonnotes} for hyperbolic subgroups $H$ of
hyperbolic groups $G$ following Gromov's suggestion about using
hyperbolic simplices (\cite{Gromov}, \S8.3). Let $X_G$ be a finite 2
complex with fundamental group $G$. Let $X_H$ be a cover of $X_G$
corresponding to the finitely presented subgroup $H$. Let $K$ be a
finite complex with fundamental group $H$. One needs to consider
homotopy equivalences between $K$ and $X_H$. Then one might try to prove
a geometric analog of Paulin's theorem \cite{Paulin} so as to obtain a
limiting action of a subgroup of $H$ on a limit metric space (in
\cite{Paulin} the limiting object is an $\mathbb{R}$-tree). This would
be an approach to answering the above question affirmatively.

The general problem attempted in this paper is one of characterizing
quasiconvexity of subgroups $H$ of hyperbolic groups $G$ purely in terms
of group theoretic notions. Swarup's question aims at one such
characterization. One might like stronger criteria, though this might be
over-optimistic. Consider the following conditions:

\begin{enumerate}
\renewcommand\labelenumi{\arabic{enumi}.}
\item $H \subset G$ is not quasiconvex.

\item $H$ has infinite height in $G$.

\item $H$ has {\em strictly infinite height} in $G$, i.e. there exist
infinitely many essentially distinct conjugates ${g_i}H{g_i^{-1}} , i =
1, 2, {\ldots}$ such that ${\cap_i}{g_i}H{g_i^{-1}} \neq \emptyset$.

\item There exists an element $g \in G$ such that $g^i \notin H$ for $i
\neq 0$ and ${\cap_i}{g^i}H{g^{-i}} \neq \emptyset$.

\item There exists an element $g \in G$ such that $g^i \notin H$ for $i
\neq 0$ and ${\cap_i}{g^i}{H_1}{g^{-i}} \neq \emptyset$ where $H_1$ is a
subgroup of $H$ isomorphic to a free product of free groups and surface
groups. 

\item There exists an element $g \in G$ such that $g^i \notin H$ for $i
\neq 0$ and ${\cap_i}{g^i}{H_1}{g^{-i}} \neq \emptyset$ where $H_1$ is a
{\em quasiconvex} subgroup of $H$ isomorphic to a free product of free
groups and surface groups.
\end{enumerate}

It is clear that $(6) \Rightarrow (5) \Rightarrow (4) \Rightarrow (3)
\Rightarrow (2) \Rightarrow (1) $ (the last implication follows from
\cite{GMRS}). One would like to know if any of these can be reversed.

\end{document}